\documentclass[12pt,reqno]{amsart}
\usepackage[latin1]{inputenc}
\usepackage[T1]{fontenc}
\usepackage{mathptmx}
\usepackage{amssymb}
\usepackage{color}

\newtheorem{theorem}{Theorem}[section]
\newtheorem{lemma}[theorem]{Lemma}
\newtheorem{corollary}[theorem]{Corollary}

\theoremstyle{definition}

\newtheorem{example}[theorem]{Example}

\numberwithin{equation}{section}

\begin{document}
\title{A note on the computation of the \\ Frobenius number of a numerical semigroup}

\author{Julio Jos\'e Moyano-Fern\'andez}

\address{Universit\"at Osnabr\"uck, FB Mathematik/Informatik, 49069
Osnabr\"uck, Germany} \email{jmoyano@uos.de}

\thanks{The author was partially supported by the Spanish Government Ministerio de Educaci\'on y Ciencia (MEC) grant MTM2007-64704 in cooperation with the European Union in the framework of the founds ``FEDER'', and by the Deutsche 
Forschungsgemeinschaft (DFG)}

\begin{abstract}
In this note we observe that the Frobenius number and therefore the conductor of a numerical semigroup can be obtained from the maximal socle degree of the quotient of the corresponding semigroup algebra by the ideal generated by the biggest generator of the semigroup.
\end{abstract}

\maketitle

\section{Introduction and Review} 

\noindent Numerical semigroups occur often in many branches of Mathematics. One of the most challenging problems in this area is the computation of the Frobenius number of the semigroup, i.e., the biggest integer not being an element of the numerical semigroup. In this paper we describe a method to calculate it based on some fundamental concepts in commutative algebra. For further details and as a general reference on numerical semigroups, the reader should refer to the works of Rosales and Garc\'{\i}a S\'{a}nchez \cite{buch}, and Ram\'{\i}rez Alfons\'{\i}n \cite{ramirez}. Much of the notation we will use comes originally from the (in many respects seminal) work of Herzog and Kunz \cite{hk}. Let $\mathbb{N}$ denote the set of nonnegative integers.
\medskip

\noindent Let $k$ be an arbitrary field. Let  $n_1, \ldots , n_d$ be positive integer numbers with $\mathrm{gcd}(n_1, \ldots , n_d)=1$. Consider the numerical semigroup 
\[
\Gamma:=\mathbb{N}n_1 + \ldots + \mathbb{N}n_d
\]
minimally generated by $n_1, \ldots , n_d$. It is well-known the existence of an element $c \in \mathbb{N}$ minimal such that $c + \mathbb{N} \subseteq \Gamma$. This number is called the \emph{conductor} of $\Gamma$, and we will denote it by $c(\Gamma)$. The number $f(\Gamma):=c(\Gamma)-1$ is then the biggest integer not belonging to $\Gamma$, and it is called the \emph{Frobenius number} of $\Gamma$.
\medskip

\noindent Let $n$ be a nonzero element of $\Gamma$. The set
\[
\mathrm{Ap}(\Gamma,n):=\{h \in \Gamma \mid h-n \notin \Gamma\}
\]
is called the \emph{Ap\'ery set} of $n$ in $\Gamma$. It is easily checked that  (cf. \cite{101})
\[
f(\Gamma)=\mathrm{max~} \mathrm{Ap}(\Gamma,n)-n. \eqno{(\dag)}
\]

\noindent Let $I \ne \varnothing$ be a subset of $\mathbb{Z}$ satisfying $I \ne \mathbb{Z}$ and $I+\Gamma \subseteq I$. Such an $I$ is said to be a \emph{fractional $\Gamma$-ideal} (sometimes also called \emph{$\Gamma$-semigroup}). The $\Gamma$-ideal $M:=\{s \in \Gamma \mid s \ne 0\}$ is the (uniquely determined) maximal ideal of $\Gamma$. It will be important in the sequel to consider also the $\Gamma$-ideal 
\[
M^{-}:=\{z \in \mathbb{Z}\mid z+M \subseteq \Gamma\}.
\]
\noindent Note that $M^{-} \supseteq \Gamma$, and since $f(\Gamma)\in M^{-}$ one has indeed $M^{-} \supsetneq \Gamma$. The inclusion $\mathbb{N} \supseteq M^{-}$ holds precisely when $\Gamma=\mathbb{N}$, and in this case $-1\in M^-$. The cardinality of the set of elements in $M^{-} \setminus \Gamma$ will be denoted by $r(\Gamma)$. Note also that 
\[
f(\Gamma) = \max \{m \mid m\in M^{-} \setminus \Gamma\}.
\]

\noindent Let $R:=k [X_1, \ldots , X_d]$ (resp. $k[t]$) be the polynomial ring over $k$ graded by $\deg (X_i)=n_i$ for every $i \in \{1, \ldots , d\}$ (resp. $\deg (t)=1$). Let $\pi$ be the graded homomorphism of $k$-algebras $\pi: R \to k[t]$ given by $X_i \mapsto t^{n_i}$ for every $i \in \{1, \ldots , d\}$. The image of $\pi$ is the semigroup ring associated with $\Gamma$, and it is denoted by $k[\Gamma]$. The homogeneous prime ideal $\mathfrak{p}:=\mathrm{ker~} \pi$ is said to be the \emph{presentation ideal} of $k[\Gamma]$ (cf. \cite{v}).
\medskip

\noindent Let us consider $\mathfrak{p}^{\prime}:=\pi_{d-1} (\mathfrak{p})$ the image in $k[X_1, \ldots , X_{d-1}]$ by the epimorphism mapping $X_d$ onto $0$, where $\pi_{d-1}$ stands for the projection onto the first $d-1$ coordinates, and define the quotient ring
\[
R^{\prime}:=k[\Gamma]/(t^{n_d}).
\]
The following ring isomorphisms are easily checked:
\[
R^{\prime} \cong k[X_1, \ldots , X_{d-1}] / \mathfrak{p}^{\prime} \cong k[\overline{X}_1, \ldots , \overline{X}_{d-1}],
\]
where $\overline{X}_i$ denotes the class of $X_i$ modulo $\mathfrak{p}^{\prime}$ for every $i \in \{1, \ldots , d-1\}$. Furthermore, the ring $R^{\prime}$ is \emph{*local}, i.e., it has a unique maximal graded ideal $\mathfrak{m}_{R^{\prime}}$.

\section{The Main Result}

\noindent Let us define the \emph{trivial submodule} (or \emph{socle}) of $R^{\prime}$ as the set of elements in $R^{\prime}$ which are annihilated by the homogeneous maximal ideal $\mathfrak{m}_{R^{\prime}}$ of $R^{\prime}$, namely
\[
\mathrm{Triv}(R^{\prime}):=\{x \in R^{\prime} \mid x \cdot \mathfrak{m}_{R^{\prime}} = (0)\}.
\]
\noindent This is the largest subspace of $R^{\prime}$ having a $R^{\prime}$-module structure of vector space, and can be identified as $\mathrm{Hom}(k,R^{\prime})$.
\medskip

\noindent Note that the set
\[
\Delta(\Gamma):=\{\lambda \in \mathbb{N} \mid \lambda + M \subseteq n_d+\Gamma\} \setminus (n_d+\Gamma),
\]
which is in fact a subset of the semigroup $\Gamma$, yields an isomorphism, say $\varphi$, between the trivial submodule $\mathrm{Triv}(R^{\prime})$ and the set of formal power series
\[
\Big \{\sum_{\nu \in \Delta(\Gamma)} r_{\nu} t^{\nu} \mid r_{\nu} \in k \Big \},
\]
whose elements are indeed polynomials. Furthermore, we have a bijection between the sets $\Delta(\Gamma)$ and $M^{-} \setminus \Gamma$ given by mapping every $\nu \in \Delta(\Gamma)$ to $\nu-n_d \in M^{-} \setminus \Gamma$. This together with the isomorphism $\varphi$ leads to the equality between the cardinality of $M^{-} \setminus \Gamma$ and the dimension of the socle:
\[
r(\Gamma)=\dim_k \mathrm{Triv}(R^{\prime}). 
\]
\noindent This means in particular that the trivial submodule $\mathrm{Triv}(R^{\prime})$ is a \emph{finite} dimensional vector space over the field $k$. Let us then choose a basis $\mathcal{B}:= \{b_1, \ldots , b_{r(\Gamma)}\}$ and take the element $\beta \in \mathcal{B}$ such that 
\[
\deg (\beta) = \max \{\deg (b_i) \mid i =1, 2, \ldots , r(\Gamma) \}.
\]

\noindent Now it is a simple matter to realise:

\begin{lemma}
The degree $\deg (\beta)$ is independent of the choice of the basis $\mathcal{B}$ of the trivial submodule of the ring $R^{\prime}$.
\end{lemma}

\noindent We are thus led to the following result:
\begin{theorem}
We have:
\[
f(\Gamma)=c(\Gamma)-1=\deg (\beta) - n_d.
\]
\end{theorem}
\begin{proof}
The proof is straightforward from the bijection $\varphi$.
\end{proof}

\begin{corollary}
We have:
\[
\mathrm{max~} \mathrm{Ap}(\Gamma,n_d)=\deg (\beta).
\]
\end{corollary}

\begin{proof}
The result follows straightforward from the equality (\dag) at the beginning of the paper. 
\end{proof}

\begin{example}
Let us take the monomial curve $C$ given by $t \mapsto (t^6,t^8,t^9)$. The corresponding numerical semigroup is $\Gamma_C=\mathbb{N} \cdot 6 +\mathbb{N} \cdot 8 +\mathbb{N} \cdot 9 $. The presentation ideal associated with $C$ is 
\[
\mathfrak{p}=(X_1^3-X_3^2, X_2^3-X_1 X_3^2),
\] 
so we have $R^{\prime}=k[X_1,X_2]/\mathfrak{p}^{\prime} \cong k [\overline{X}_1, \overline{X}_2]$ with $\mathfrak{p}^{\prime} = (X_1^3,X_2^3)$ and $\overline{X}_i \equiv X_i \mathrm{~mod~} \mathfrak{p}^{\prime}$ for $i=1,2$. Therefore we get $\mathrm{Triv}(R^{\prime})= k \cdot \overline{X}_1^2 \overline{X}_2^2$, and it is clearly seen that $f(\Gamma_C)=\deg (\overline{X}_1^2 \overline{X}_2^2)-9=12+16-9=19$; hence $c(\Gamma_C)=20$, as one might have also checked---easily in this example---straight from the semigroup $\Gamma_C$.
\end{example}

\end{document}